\documentclass{article}
\usepackage{amssymb}
\usepackage{amsmath}

\begin{document}
%
%
\title{On Systems of Linear Quaternion Functions}
\author{Todd A. Ell}

\maketitle      

\begin{abstract}
A method of reducing general quaternion functions of first degree, i.e., linear quaternion functions, to \emph{quaternary canonical form} is given. Linear quaternion functions, once reduced to canonical form, can be maintained in this form under functional composition. Furthermore, the composition operation is symbolically identical to quaternion multiplication, making manipulation and reduction of systems of linear quaternion functions straight forward.
\end{abstract}    

%
\section{Introduction}
\label{Intro}
Real linear functions take the monomial form: $f(x)~=~mx$, where $x,m \in\mathbb{R}$. By contrast the general quaternion linear function has a multi-nomial form:
\begin{equation}
\label{QuatLinear}
f\left( q \right) = \sum\limits_{p = 1}^P {m_p \,q\,n_p } ,
\end{equation}
where all factors are quaternion valued, i.e., $q, m_p ,n_p \in \mathbb{H}$. The number of terms in the summation can be arbitrarily large due to quaternion multiplication being non-commutative.    

Unlike Hardy \cite[p.199-200]{Hardy:1} who studied the inversion of these functions or Hathaway \cite[chp. 4]{Hathaway:1} who was interested in their properties, this work addresses the reduction of \emph{systems} or networks of such functions, formed by scalar-weighted sums and functional composition.  The number of monomial terms in such systems quickly become intractable. 

Herein is a method for reducing such systems in an efficient and transparent manner.  Two reduction methods are key to achieving this goal.  First is conversion of general linear quaternion functions, as given in equation (\ref{QuatLinear}), to at most a four-tuple form.  It will be shown that equation (\ref{QuatLinear}) can always be reduced to the form:
\begin{equation}
\label{LQS}
f\left( q \right) = A\,q + B\,q\,i + C\,q\,j + D\,q\,k,
\end{equation}
where $A,B,C,D \in \mathbb{H}$ and $i,j,k$ are the standard hypercomplex operators. Linear quaternion functions of this type shall be said to be in \emph{quaternary canonical form}. Just as a quaternion can be associated with a 4-tuple of \emph{reals} as $a1 + bi + cj + dk \Leftrightarrow \left( {a,b,c,d} \right)$, likewise the canonical linear function can be associated with a 4-tuple of quaternions as
\[
A\,q + B\,q\,i + C\,q\,j + D\,q\,k \Leftrightarrow \left\{ {A,B,C,D} \right\}.
\]
The four-tuple of quaternions will be used as a shorthand notation for the canonical form.

The linear sum of two such functions, say $f_1$ and $f_2$, is immediately reducible to canonical form, i.e., 
\[
f_3 \left( q \right) = f_1 \left( q \right) + f_2 \left( q \right) = A_3 q + B_3 q\,i + C_3 q\,j + D_3 q\,k
\]
where 
\[
A_3  = A_1  + A_2 ,\quad B_3  = B_1  + B_2 ,\quad C_3  = C_1  + C_2 ,\quad D_3  = D_1  + D_2 .
\]

However, functional composition of two such functions naturally expand into sixteen terms, which require reduction back to four terms. The number of terms in multiple compositions expand exponentially.  The second reduction method provides a simple composition rule which automatically yields a quaternary form, provided the two functions being composed are already in this form.    

The next section provides preliminary material to introduce notation and nomenclature used in this work, as well as, matrix  and involution operations useful in deriving the main results. Section \ref{ReductionToCannonicalForm} presents the method of reduction to quaternary form using either equivalent matrix operations or simpler involution operators.  The matrix-based version provides a numerically efficient algorithm, while the involution-based version is better suited to symbolic analysis.  Section \ref{SystemReduction} presents the composition algorithm. 


Motivation behind this work is in the use of arrays of linear quaternion functions as filtering operations on quaternion valued data arrays. This allows for complex geometric operations on three- and four-dimensional data sets.  Only after the basic manipulation of systems of such functions are detailed can the more difficult problem of linear quaternion filters be addressed.
 
%
\section{Preliminaries}
\label{Preliminaries}
\subsection{Quaternions}
The \emph{quaternions}, $\mathbb{H}$, are a four dimensional, non-commutative division algebra over the reals, $\mathbb{R}$. Multiplication is defined by the hypercomplex operator rules
\[
i^2  = j^2  = k^2  = ijk =  - 1.
\]

The typical quaternion is written in Cartesian form as 
\[
p = p_0  + p_1 i + p_2 j + p_3 k
\]
where $p_{0,1,2,3} \in \mathbb{R}$. 

The quaternion \emph{conjugate}, $\bar p$, is defined as negating the non-real part, i.e., 
\[
\bar p = p_0  - p_1 i - p_2 j - p_3 k
\]

\subsection{Quaternion-Matrix Homomorphism }
\label{Quaternion.matrix.Homomorphism}
The equivalence between quaternion and 4$\times$4 real matrix multiplication was first established by Peirce \cite{Peirce:1}, later verified by Spottiswoode \cite{Spottiswoode:1}, made explicit by Cayley \cite{Cayley:1,Cayley:2} and later in Tait's book \cite{Tait:1}. Walker \cite{Walker:1} provides a concise introduction and Ickes \cite{Ickes:1} uses this equivalence to provide fast numerical computations on linear quaternion equations. 

These quaternion matrix representations are not unique. Farebrother et al. \cite{Farebrother:1} identified 48 different representations, not considering all automorphisms of the quaternion algebra.  Later, Ducati \cite{Ducati:1} divided these representations into two groups based on the left and right action of the hypercomplex operators.  Mackey \cite{Mackey:1} used this concept to solve the eigenvalue problem.  Fassbender et al. \cite{Fassbender:1} used them to solve the Hamiltonian eigenproblem.   DeLeo and Rotelli \cite{DeLeo:1} used the same concept applied to complex 2$\times$2 matrices.   

Using the hypercomplex operator product rules, the quaternion product of two arbitrary quaternions, $q,p \in \mathbb{H}$, written in Cartesian form as 
\[
\begin{gathered}
  p = p_0  + p_1 i + p_2 j + p_3 k \hfill \\
  q = q_0  + q_1 i + q_2 j + q_3 k \hfill \\ 
\end{gathered} 
\]
is 
\[
\begin{gathered}
  r_0  = p_0 q_0  - p_1 q_1  - p_2 q_2  - p_3 q_3  \hfill \\
  r_1  = p_1 q_0  + p_0 q_1  - p_3 q_2  + p_2 q_3  \hfill \\
  r_2  = p_2 q_0  + p_0 q_2  + p_3 q_1  - p_1 q_3  \hfill \\
  r_3  = p_3 q_0  + p_0 q_3  - p_2 q_1  + p_1 q_2  \hfill \\ 
\end{gathered} 
\]

where $r = pq = r_0 + r_1i + r_2j + r_3k$. Gathering terms into matrix-vector notation we obtain two distinct matrix-vector products: the \emph{standard} form
\begin{equation}
\label{R=[p]Q}
\left[ {\begin{array}{*{20}c}
   {r_0 }  \\
   {r_1 }  \\
   {r_2 }  \\
   {r_3 }  \\

 \end{array} } \right] = \left[ {\begin{array}{*{20}c}
   {p_0 } & { - p_1 } & { - p_2 } & { - p_3 }  \\
   {p_1 } & {p_0 } & { - p_3 } & {p_2 }  \\
   {p_2 } & {p_3 } & {{\text{ }}p_0 } & { - p_1 }  \\
   {p_3 } & { - p_2 } & {p_1 } & {p_0 }  \\

 \end{array} } \right]\left[ {\begin{array}{*{20}c}
   {q_0 }  \\
   {q_1 }  \\
   {q_2 }  \\
   {q_3 }  \\

 \end{array} } \right]
\end{equation}
which preserves the order of the product, and the \emph{transmuted} form
\begin{equation}
\label{R=[q]P}
\left[ {\begin{array}{*{20}c}
   {r_0 }  \\
   {r_1 }  \\
   {r_2 }  \\
   {r_3 }  \\

 \end{array} } \right] = \left[ {\begin{array}{*{20}c}
   {q_0 } & { - q_1 } & { - q_2 } & { - q_3 }  \\
   {q_1 } & {q_0 } & {q_3 } & { - q_2 }  \\
   {q_2 } & { - q_3 } & {q_0 } & {q_1 }  \\
   {q_3 } & {q_2 } & { - q_1 } & {q_0 }  \\

 \end{array} } \right]\left[ {\begin{array}{*{20}c}
   {p_0 }  \\
   {p_1 }  \\
   {p_2 }  \\
   {p_3 }  \\

 \end{array} } \right]
\end{equation}
which reverses the order. 

Following Ickes \cite{Ickes:1}, let the standard matrix form be denoted $\left[ p\right ]$ and the transmuted matrix form as $\left[ q \right]^\dag$, where the `$\dag$' signifies transposing the lower-right 3$\times$3 sub-matrix.  

In the quaternion product $r=pq$, if we regard $p$ as a left-handed operator on the state variable $q$, then the $p$-operator is encoded as a 4$\times$4 real matrix and the $q$-state as a real $4\times$1 column vector.  Conversely, if $q$ is regarded as a right-handed operator on the state variable $p$, then the operator is again encoded as a real matrix (transmuted to the left) and the state as a real vector.  

\subsection{Quaternion Operators}
\label{BarOp}

Since left and right quaternion multiplications are unique, we can view various combinations as independent operators.  Following De Leo and Rotelli \cite{DeLeo:1}, one can define a \emph{bar}-operator notation to provide a shorthand representation for these combinations. The \emph{bar}-operator, denoted $e_1 | e_2$, is defined as:
\[
e_1 |e_2 \left( q \right) \to e_1 \left( q \right)e_2 ,{\text{     where }}e_n  \in \left\{ {1,i,j,k} \right\}
\]
so there are sixteen such operators as listed in tabular form below:
\[
\begin{array}{*{20}c}
   {1|1,} & {1|i,} & {1|j,} & {1|k,}  \\
   {i|1,} & {i|i,} & {i|j,} & {i|k,}  \\
   {j|1,} & {j|i,} & {j|j,} & {j|k,}  \\
   {k|1,} & {k|i,} & {k|j,} & {k|k.}  \\

 \end{array} 
\]

These quaternion operators are consistent with the standard hypercomplex operators and have similar operator formulae
\[
\begin{gathered}
  \left( {i|1} \right)^2  = \left( {j|1} \right)^2  = \left( {k|1} \right)^2  = \left( {i|1} \right)\left( {j|1} \right)\left( {k|1} \right) =  - 1|1 \hfill \\
  \left( {1|i} \right)^2  = \left( {1|j} \right)^2  = \left( {1|k} \right)^2  = \left( {1|i} \right)\left( {1|j} \right)\left( {1|k} \right) =  - 1|1 \hfill \\ 
\end{gathered} 
\]
so that various operator equations can be manipulated directly, e.g., $\left( {i|j} \right)\left( {j|k} \right) = \left( {ij|kj} \right) =  - \left( {k|i} \right)$.  

Using the standard and transmuted matrix homomorphisms given in the previous section, i.e., $e_1 |e_2  \to \left[ {e_1 } \right]\left[ {e_2 } \right]^\dag $, each can be written as a matrix. For example, 
\[
\left( {k|i} \right) \to \left[ k \right]\left[ i \right]^\dag   = \left[ {\begin{array}{*{20}c}
   0 & 0 & { + 1} & 0  \\
   0 & 0 & 0 & { - 1}  \\
   { + 1} & 0 & 0 & 0  \\
   0 & { - 1} & 0 & 0  \\

 \end{array} } \right].
\]

The basic quaternion operator formulae can be extended by forming weighted linear combinations. Linearity of the quaternion operations implies that for four arbitrary quaternions, $A,B,C,D\in \mathbb{H}$, written in component form as
\[
\begin{gathered}
  A = \alpha _0  + \alpha _1 i + \alpha _2 j + \alpha _3 k, \hfill \\
  B = \beta _0  + \beta _1 i + \beta _2 j + \beta _3 k, \hfill \\
  C = \gamma _0  + \gamma _1 i + \gamma _2 j + \gamma _3 k, \hfill \\
  D = \delta _0  + \delta _1 i + \delta _2 j + \delta _3 k, \hfill \\ 
\end{gathered} 
\]
when combined with four operators $1|1$, $1|i$, $1|j$, and $1|k$, respectively, form the four extended operators 
\[
\begin{gathered}
  \left( {A|1} \right) = \alpha _0 \left( {1|1} \right) + \alpha _1 \left( {i|1} \right) + \alpha _2 \left( {j|1} \right) + \alpha _3 \left( {k|1} \right), \hfill \\
  \left( {B|i} \right) = \beta _0 \left( {1|i} \right) + \beta _1 \left( {i|i} \right) + \beta _2 \left( {j|i} \right) + \beta _3 \left( {k|i} \right), \hfill \\
  \left( {C|j} \right) = \gamma _0 \left( {1|j} \right) + \gamma _1 \left( {i|j} \right) + \gamma _2 \left( {j|j} \right) + \gamma _3 \left( {k|j} \right), \hfill \\
  \left( {D|k} \right) = \delta _0 \left( {1|k} \right) + \delta _1 \left( {i|k} \right) + \delta _2 \left( {j|k} \right) + \delta _3 \left( {k|k} \right). \hfill \\ 
\end{gathered} 
\]

In the standard case the quaternions $\pm 1$, $\pm i$, $\pm j$, and $\pm k$ form a non-abelian group of order eight, with multiplication as the group operator.  In the quaternion-operator case, these become $\pm 1|1$, $\pm 1|i$, $\pm 1|j$, $\pm 1|k$, $\pm i|1$, $\pm j|1$, and $\pm k|1$ giving a 14th order group with composition being the group operator.

\subsection{Quaternion Involutions }
\label{Quaternion.Involutions}
A function, $f$, is an \emph{involution}\footnote{This is also known as an \emph{involutory anti-automorphism}.} 
 if it has the following properties \cite{Hazewinkel:1}: 
\begin{enumerate}
 \item is it's own inverse: $f(f(q)) = q$,  
 \item is linear:  $f(q_1+q_2) = f(q_1) + f(q_2)$, and 
 \item is anti-homomorphic: $f(q_1q_2) = f(q_2) f(q_1)$.
\end{enumerate} 

A function is an \emph{anti-involution}\footnote{Likewise, this is also known as an \emph{involutory automorphism}.} if the third property, above, is homomorphic, i.e.,  $f(q_1q_2) = f(q_1) f(q_2)$.  

Using the quaternion bar-operators in conjunction with conjugation provides a \emph{generalized conjugate},  
\[
\bar q^\varepsilon   \to  - \varepsilon \bar q\varepsilon ,\quad \varepsilon  \in \left\{ {i,j,k} \right\}
\]
for any $q \in \mathbb{H}$, which are involutions under these definitions \cite{Benn:1}.  The conjugate of these functions,   
\[
\overline {\bar q^\varepsilon  }  \to  - \varepsilon q\varepsilon ,\quad \varepsilon  \in \left\{ {i,j,k} \right\}
\]
are anti-involutions \cite{Ell:1}.

For an arbitrary quaternion $q=q_0+iq_1+jq_2+kq_3$ the following identities hold in terms of involutions
\begin{equation}
\label{Inv}
\begin{gathered}
  q_1  = \tfrac{1}
{2}\left( {\bar q^i  + \bar q} \right)i \\ 
  q_2  = \tfrac{1}
{2}\left( {\bar q^j  + \bar q} \right)j \\ 
  q_3  = \tfrac{1}
{2}\left( {\bar q^k  + \bar q} \right)k \\ 
\end{gathered} 
\end{equation}
which are analogous to the formula for the scalar part of a quaternion using the standard conjugate
\begin{equation}
\label{Sca}
q_0  = \tfrac{1}
{2}\left( {q + \bar q} \right)
\end{equation}
except they extract the remaining real-valued component-parts of any quaternion.   

Finally, it should be pointed out that all of these functions (generalized conjugates, and component-part functions) are readily placed in quaternary form. This is achieved by first expressing the standard quaternion conjugate in involutory form, i.e.,
\[
\bar q = \tfrac{1}
{2}\left( { - q + iqi + jqj + kqk} \right).
\]

Hence its canonical four-tuple is $\bar q = \left\{ { - \tfrac{1}{2},\tfrac{i}{2},\tfrac{j}{2},\tfrac{k}{2}} \right\}$. Proof of this involutory form is given in \cite{Ell:1}.

%
\section{Reduction to Canonical Form}
\label{ReductionToCannonicalForm}
\subsection{Matrix Method}
At its core, the matrix method is the process for taking an arbitrary operator matrix and converting it to the sum of four operator matrices, each matrix equivalent to an operator quaternion.  

The matrix-quaternion homomorphism of section \ref{Quaternion.matrix.Homomorphism} suggest that the individual terms in equation (\ref{QuatLinear}) can be converted into matrix form as
\[
m_p \left( {qn_p } \right) \to \left[ {m_p } \right]\left( {\left[ {n_p } \right]^\dag  \vec q} \right) = \left[ {m_p } \right]\left[ {n_p } \right]^\dag  \vec q
\]
where $\vec q$ denotes the state-vector form of quaternion $q$.

If a numerical reduction to the problem is desired, then the following matrix equation could be used
\[
f\left( q \right) \to \left\{ {\sum\limits_{p = 1}^P {\left[ {m_p } \right]\left[ {n_p } \right]^\dag  } } \right\}\vec q.
\]

Even though the resulting 4-vector represents a quaternion, the intermediate matrices $\left[ {m_p } \right]\left[ {n_p } \right]^\dag$ cannot be decoded into a equivalent \emph{single} quaternion since in general it does not have the correct pattern of elements. A general 4$\times$4 real matrix has sixteen degrees of freedom, whereas a quaternion has only four (whether or not it is in matrix form). So arbitrary real matrices cannot be written as a \emph{single} quaternion. However, it can be written as the direct sum of four matrices, each matrix equivalent to a quaternion. 

Let $A,B,C,D \in \mathbb{H}$, be defined as
\[
\begin{gathered}
  A = \alpha _0  + \alpha _1 i + \alpha _2 j + \alpha _3 k, \hfill \\
  B = \beta _0  + \beta _1 i + \beta _2 j + \beta _3 k, \hfill \\
  C = \gamma _0  + \gamma _1 i + \gamma _2 j + \gamma _3 k, \hfill \\
  D = \delta _0  + \delta _1 i + \delta _2 j + \delta _3 k, \hfill \\ 
\end{gathered} 
\]
then each of the above quaternions contains four degrees of freedom. Hence the combined equation  
\[
\left( {A|1} \right) + \left( {B|i} \right) + \left( {C|j} \right) + \left( {D|k} \right)
\]
has a total of sixteen; the number necessary to encode an arbitrary 4$\times$4 matrix. Writing each quaternion operator in matrix form as
\[
\left( {A|1} \right) \to \left[ A \right]\left[ 1 \right]^\dag   = \left[ {\begin{array}{*{20}c}
   {\alpha _0 } & { - \alpha _1 } & { - \alpha _2 } & { - \alpha _3 }  \\
   {\alpha _1 } & {\alpha _0 } & { - \alpha _3 } & {\alpha _2 }  \\
   {\alpha _2 } & {\alpha _3 } & {\alpha _0 } & { - \alpha _1 }  \\
   {\alpha _3 } & { - \alpha _2 } & {\alpha _1 } & {\alpha _0 }  \\

 \end{array} } \right],
\]
\[
\left( {B|i} \right) \to \left[ B \right]\left[ i \right]^\dag   = \left[ {\begin{array}{*{20}c}
   { - \beta _1 } & { - \beta _0 } & {\beta _3 } & { - \beta _2 }  \\
   {\beta _0 } & { - \beta _1 } & { - \beta _2 } & { - \beta _3 }  \\
   {\beta _3 } & { - \beta _2 } & {\beta _1 } & {\beta _0 }  \\
   { - \beta _2 } & { - \beta _3 } & { - \beta _0 } & {\beta _1 }  \\

 \end{array} } \right],
\]
\[
\left( {C|j} \right) \to \left[ C \right]\left[ j \right]^\dag   = \left[ {\begin{array}{*{20}c}
   { - \gamma _2 } & { - \gamma _3 } & { - \gamma _0 } & {\gamma _1 }  \\
   { - \gamma _3 } & {\gamma _2 } & { - \gamma _1 } & { - \gamma _0 }  \\
   {\gamma _0 } & { - \gamma _1 } & { - \gamma _2 } & { - \gamma _3 }  \\
   {\gamma _1 } & {\gamma _0 } & { - \gamma _3 } & {\gamma _2 }  \\

 \end{array} } \right],
\]
and
\[
\left( {D|k} \right) \to \left[ D \right]\left[ k \right]^\dag   = \left[ {\begin{array}{*{20}c}
   { - \delta _3 } & {\delta _2 } & { - \delta _1 } & { - \delta _0 }  \\
   {\delta _2 } & {\delta _3 } & {\delta _0 } & { - \delta _1 }  \\
   { - \delta _1 } & { - \delta _0 } & {\delta _3 } & { - \delta _2 }  \\
   {\delta _0 } & { - \delta _1 } & { - \delta _2 } & { - \delta _3 }  \\

 \end{array} } \right],
\]
then the combined operator in matrix form is
\[
\begin{gathered}
  \left( {A|1} \right) + \left( {B|i} \right) + \left( {C|j} \right) + \left( {D|k} \right) \to  \hfill \\
  \quad \left[ {\begin{array}{*{20}c}
   {\alpha _0  - \beta _1  - \gamma _2  - \delta _3 } & {\;\; - \alpha _1  - \beta _0  - \gamma _3  + \delta _2 }  \\
   {\alpha _1  + \beta _0  - \gamma _3  + \delta _2 } & {\;\;\alpha _0  - \beta _1  + \gamma _2  + \delta _3 }  \\
   {\alpha _2  + \beta _3  + \gamma _0  - \delta _1 } & {\;\;\alpha _3  - \beta _2  - \gamma _1  - \delta _0 }  \\
   {\alpha _3  - \beta _2  + \gamma _1  + \delta _0 } & {\;\; - \alpha _2  - \beta _3  + \gamma _0  - \delta _1 }  \\

 \end{array} } \right. \cdots  \hfill \\
  \quad \quad \quad \quad \left. {\begin{array}{*{20}c}
   { - \alpha _2  + \beta _3  - \gamma _0  - \delta _1 } & {\;\; - \alpha _3  - \beta _2  + \gamma _1  - \delta _0 }  \\
   { - \alpha _3  - \beta _2  - \gamma _1  + \delta _0 } & {\;\;\alpha _2  - \beta _3  - \gamma _0  - \delta _1 }  \\
   {\alpha _0  + \beta _1  - \gamma _2  + \delta _3 } & {\;\; - \alpha _1  + \beta _0  - \gamma _3  - \delta _2 }  \\
   {\alpha _1  - \beta _0  - \gamma _3  - \delta _2 } & {\;\;\alpha _0  + \beta _1  + \gamma _2  - \delta _3 }  \\

 \end{array} } \right]. \hfill \\ 
\end{gathered} 
\]

Setting this equal to an arbitrary matrix, $R$, with elements $\left\{ {r_{ij} } \right\}$, and solving for the components of the quaternions one finds
\begin{equation}
\label{A1}
A = \left[ {\begin{array}{*{20}c}
   {\alpha _0 }  \\
   {\alpha _1 }  \\
   {\alpha _2 }  \\
   {\alpha _3 }  \\

 \end{array} } \right] = \frac{1}
{4}\left[ {\begin{array}{*{20}c}
   { + r_{11}  + r_{22}  + r_{33}  + r_{44} }  \\
   { + r_{21}  - r_{12}  + r_{43}  - r_{34} }  \\
   { + r_{31}  - r_{13}  - r_{42}  + r_{24} }  \\
   { + r_{41}  - r_{14}  + r_{32}  - r_{23} }  \\

 \end{array} } \right],
\end{equation}

\begin{equation}
\label{Bi}
B = \left[ {\begin{array}{*{20}c}
   {\beta _0 }  \\
   {\beta _1 }  \\
   {\beta _2 }  \\
   {\beta _3 }  \\

 \end{array} } \right] = \frac{1}
{4}\left[ {\begin{array}{*{20}c}
   { + r_{21}  - r_{12}  - r_{43}  - r_{34} }  \\
   { - r_{11}  - r_{22}  + r_{33}  + r_{44} }  \\
   { - r_{41}  - r_{14}  - r_{32}  - r_{23} }  \\
   { + r_{31}  + r_{13}  - r_{42}  - r_{24} }  \\

 \end{array} } \right],
\end{equation}

\begin{equation}
\label{Cj}
C = \left[ {\begin{array}{*{20}c}
   {\gamma _0 }  \\
   {\gamma _1 }  \\
   {\gamma _2 }  \\
   {\gamma _3 }  \\

 \end{array} } \right] = \frac{1}
{4}\left[ {\begin{array}{*{20}c}
   { + r_{31}  - r_{13}  + r_{42}  - r_{24} }  \\
   { + r_{41}  + r_{14}  - r_{32}  - r_{23} }  \\
   { - r_{11}  + r_{22}  - r_{33}  + r_{44} }  \\
   { - r_{21}  - r_{12}  - r_{43}  - r_{34} }  \\

 \end{array} } \right],
\end{equation}

and
\begin{equation}
\label{Dk}
D = \left[ {\begin{array}{*{20}c}
   {\delta _0 }  \\
   {\delta _1 }  \\
   {\delta _2 }  \\
   {\delta _3 }  \\

 \end{array} } \right] = \frac{1}
{4}\left[ {\begin{array}{*{20}c}
   { + r_{41}  - r_{14}  - r_{32}  + r_{23} }  \\
   { - r_{31}  - r_{13}  - r_{42}  - r_{24} }  \\
   { + r_{21}  + r_{12}  - r_{43}  - r_{34} }  \\
   { - r_{11}  + r_{22}  + r_{33}  - r_{44} }  \\

 \end{array} } \right].
\end{equation}

This is not as difficult as it appears since the problem decomposes into four sets of four unknowns with each set independent of the others, e.g., the coefficient set $\left\{ {\alpha _0 ,\beta _1 ,\gamma _2 ,\delta _3 } \right\}$ is solved using only the diagonal terms in both matrices. 

The steps in reducing equation (\ref{QuatLinear}) to quaternary form are:

\begin{enumerate}
 \item Map coefficients $m_p$ and $n_p$ to matrix-space using equations (\ref{R=[p]Q}-\ref{R=[q]P})
\[
m_p  \to \left[ {m_p } \right],\quad n_p  \to \left[ {n_p } \right]^\dag  .
\]
\item Construct single matrix $R$ as
\[
{\mathbf{R}} = \sum\limits_{p = 1}^P {\left[ {m_p } \right]\left[ {n_p } \right]^\dag  } .
\]
\item Construct equivalent quaternion factors using equations (\ref{A1}-\ref{Dk})
\[
R \to \left( {A|1} \right) + \left( {B|i} \right) + \left( {C|j} \right) + \left( {D|k} \right).
\]
\end{enumerate}

\subsection{Involution Method}
The matrix equivalent method outlined above is useful for numerical computations, however, it requires manipulation at the component level of each quaternion, $m_p$ and $n_p$. The involution method, outlined in this section, provides an alternative solution in terms of $m_p$ and $n_p$ directly, without these explicit references, making it more suitable for symbolic analysis.  

Starting with equation (\ref{QuatLinear}), let each $n_p$ be expanded into component form as
\[
n_p  = w_p  + ix_p  + jy_p  + kz_p 
\]
then the full equation can be expanded as 
\[
\begin{gathered}
  f\left( q \right) = \sum\limits_{p = 1}^P {m_p \,q\,\left( {w_p  + ix_p  + jy_p  + kz_p } \right)}  \hfill \\
  \quad \quad  = \sum\limits_{p = 1}^P {m_p w_p q}  + \sum\limits_{p = 1}^P {m_p x_p q\,i}  \hfill \\
  \quad \quad \quad \quad \quad \quad  + \sum\limits_{p = 1}^P {m_p y_p q\,j}  + \sum\limits_{p = 1}^P {m_p z_p q\,k}  \hfill \\
  \quad \quad  = \left[ {\sum\limits_{p = 1}^P {m_p w_p } } \right]q + \left[ {\sum\limits_{p = 1}^P {m_p x_p } } \right]q\,i \hfill \\
  \quad \quad \quad \quad  + \left[ {\sum\limits_{p = 1}^P {m_p y_p } } \right]q\,j + \left[ {\sum\limits_{p = 1}^P {m_p z_p } } \right]q\,k. \hfill \\ 
\end{gathered} 
\]
Now let
\[
\begin{gathered}
  A = \sum\limits_{p = 1}^P {m_p w_p } ,\quad B = \sum\limits_{p = 1}^P {m_p x_p } , \hfill \\
  C = \sum\limits_{p = 1}^P {m_p y_p } ,\quad D = \sum\limits_{p = 1}^P {m_p z_p } , \hfill \\ 
\end{gathered} 
\]
so that 
\[
f\left( q \right) = Aq + Bq\,i + Cq\,j + Dq\,k.
\]
The problem with this formula is that we are required to decompose quaternions $n_p$ into component form, whereas we are looking for a formula that does not require this step. The answer to this conundrum is to apply the involution formulas (\ref{Inv}-\ref{Sca}).  Namely, let
\[
\begin{gathered}
  w_p  = \tfrac{1}
{2}\left( {n_p  + \bar n_p } \right),\;\;\quad x_p  = \tfrac{1}
{2}\left( {\bar n_p^i  + \bar n_p } \right)i, \hfill \\
  y_p  = \tfrac{1}
{2}\left( {\bar n_p^j  + \bar n_p } \right)j,\quad z_p  = \tfrac{1}
{2}\left( {\bar n_p^k  + \bar n_p } \right)k, \hfill \\ 
\end{gathered} 
\]
so that 
\[
\begin{gathered}
  A = \tfrac{1}
{2}\sum\limits_{p = 1}^P {m_p \left( {n_p  + \bar n_p } \right)} , \hfill \\
  B = \tfrac{1}
{2}\sum\limits_{p = 1}^P {m_p \left( {\bar n_p^i  + \bar n_p } \right)i} , \hfill \\
  C = \tfrac{1}
{2}\sum\limits_{p = 1}^P {m_p \left( {\bar n_p^j  + \bar n_p } \right)j} , \hfill \\
  D = \tfrac{1}
{2}\sum\limits_{p = 1}^P {m_p \left( {\bar n_p^k  + \bar n_p } \right)k} . \hfill \\ 
\end{gathered} 
\]
 
%
\section{System Reduction}
\label{SystemReduction}
One could use the methods of the previous section to reduce compositions of two linear functions, i.e., where $f_2 \left( q \right) \circ f_1 \left( q \right) = f_2 \left( {f_1 \left( q \right)} \right)$. What is not apparent is that functional composition can be reduced in a more elegant way.  This alternative lends itself to better insight into what happens in these composition operations.  

Using bar-operator notation, composition is reduced to ordered multiplication, hence one can build a {\it composition} table showing how components of the composed functions are reduced.  The composition table for $f_2 \circ f_1$ is given below
\[
\begin{array}{*{20}c}
   {f_2  \circ f_1 } &\vline &  {A_1 |1} & {B_1 |i} & {C_1 |j} & {D_1 |k}  \\
\hline
   {A_2 |1} &\vline &  { + A_2 A_1 |1} & { + A_2 B_1 |i} & { + A_2 C_1 |j} & { + A_2 D_1 |k}  \\
   {B_2 |i} &\vline &  { + B_2 A_1 |i} & { - B_2 B_1 |1} & { - B_2 C_1 |k} & { + B_2 D_1 |j}  \\
   {C_2 |j} &\vline &  { + C_2 A_1 |j} & { + C_2 B_1 |k} & { - C_2 C_1 |1} & { - C_2 D_1 |i}  \\
   {D_2 |k} &\vline &  { + D_2 A_1 |k} & { - D_2 B_1 |j} & { + D_2 C_1 |i} & { - D_2 D_1 |1}  \\

 \end{array} 
\]

Collecting like terms, we find that 
\[
f_2  \circ f_1 \left( q \right) = A_3 q + B_3 q\,i + C_3 q\,j + D_3 q\,k
\]
where 
\[
\begin{gathered}
  A_3  = A_2 A_1  - B_2 B_1  - C_2 C_1  - D_2 D_1 , \hfill \\
  B_3  = A_2 B_1  + B_2 A_1  - C_2 D_1  + D_2 C_1 , \hfill \\
  C_3  = A_2 C_1  + B_2 D_1  + C_2 A_1  - D_2 B_1 , \hfill \\
  D_3  = A_2 D_1  - B_2 C_1  + C_2 B_1  + D_2 A_1 . \hfill \\ 
\end{gathered} 
\]

Careful examination of the 4-tuple resulting from a product of two quaternions with the 4-tuple resulting from the functional composition of two canonical equations reveals that they are transmuted versions of the same operational formula. This is a consequence of the bar-operators having a similar operator formula as the standard hypercomplex operator formula.  If the canonical form had been defined with the full quaternions on the right, e.g.,
\[
f\left( q \right) = q\,A' + i\,q\,B' + j\,q\,C' + k\,q\,D',
\]
then the composition formula would have been identical to the quaternion multiplication formula.

\section{Conclusions}%
The goal of this work is the reduction and manipulation of systems of linear quaternion functions. To this end two results have been provided.  First, even though a general linear quaternion function may contain an arbitrary number of terms, it can always be reduced to at most four terms.  This reduced function, called the \emph{quaternary canonical form} of the function, is completely specified with a four-tuple of quaternions.  Second, series (cascaded composition) and parallel (weighted-sum) combinations of such canonical functions can be written directly in canonical form without recourse to the previous reduction method.  


\end{document}